\documentclass[a4paper,12pt]{amsart}
\usepackage[T1]{fontenc}
\usepackage{amssymb}
\usepackage{amsopn} 
\usepackage[square]{harvard}
\usepackage{amsthm}
\usepackage{amsmath}

\newcommand{\N}{{\mathbb N}}
\newcommand{\Z}{{\mathbb Z}}
\newcommand{\Q}{{\mathbb Q}}
\newcommand{\R}{{\mathbb R}}
\newcommand{\C}{{\mathbb C}}
\let \op = \operatorname
\newcommand{\SL}{\operatorname{SL_2}}
\newcommand{\spn}{\operatorname{span}}

\newtheorem{Lemma}{Lemma}[section]
\newtheorem{Prop}{Proposition}[section]

\newtheorem{Theorem}{Theorem}[section]
\newtheorem{Corollary}{Corollary}[section]

\theoremstyle{remark}
\newtheorem{Remark}{Remark}[section]
\newtheorem{Example}{Example}[section]

\begin{document}
\title{Asymptotic densities of newforms}
\author{Morten Skarsholm Risager}
\date{November 13, 2001}
\begin{abstract}We define the counting function for non-analytic (Maass) newforms of Hecke
  congruence groups $\Gamma_{\!0}(M)$. We then calculate the three main terms of this
  counting function and give necessary and sufficient conditions on
  $M$ for this counting function to have the same shape as if it were
  counting eigenvalues related to a cocompact group.  
\end{abstract}
\address{Department of Mathematical Sciences, University of Aarhus
Ny Munkegade Building 530, 8000 Århus, Denmark
}
\email{risager@imf.au.dk}

\maketitle
\section{Introduction}
Let $\Gamma$ be a congruence subgroup of the full modular group. It is
well known that the selfadjoint automorphic Laplacian, $\Delta_\Gamma$,
has infinitely many eigenvalues, 
$$0=\lambda_0\leq\lambda_1^\Gamma\leq\ldots\leq\lambda_i^\Gamma\leq\ldots,$$
listed with their multiplicities which are finite. Selberg has proved that the counting
function $$N_\Gamma(\lambda)= \# \{i|\lambda_i^\Gamma\leq\lambda\}$$
satisfies a Weyl law namely
\begin{equation}\label{weyl} N_\Gamma(\lambda)=\frac{|F_\Gamma|}{4\pi}\lambda+O(\sqrt{\lambda}\log \lambda),\end{equation}
where $|F_\Gamma|$ is the area of a fundamental domain of
$\Gamma$. There are various refinements of this result (see
e.g. theorem \ref{refined}).

In this paper we investigate what happens if we only count the
eigenvalues corresponding to newforms.
\section{Newforms and oldforms}
The theory of newforms was originally developed by
\citeasnoun{MR42:3022} for holomorphic forms. Their theory can be
translated into a similar theory of Maass forms which are the ones we
are studying. This has been done independently by various people and details may be
found in e.g. \cite{andreas}. We shall only need one result (Lemma
\ref{rekurrens} below) and shall hence only sketch enough of the
theory for this result to make sense.

For any $\lambda>0$ , $M\in\N$ we denote by  $A(\lambda,M)$ the
 $\lambda$-eigenspace for $\Delta_{\Gamma_{\!0}(M)}$, where
 $\Gamma_{\!0}(M)$ is the Hecke congruence group of level $M$ i.e. 
 $$\Gamma_{\!0}(M)=\left\{\gamma\in \SL(\Z)\left|\gamma=\left(\begin{array}{cc}a&b\\c&d\end{array}\right) \right. c\equiv 0 \mod
     M\right\}.$$ 

Then it is obvious that
 \begin{equation}\label{omskriv}N_{\Gamma_{\!0}(M)}(\lambda)=1+\sum_{0<\tilde\lambda\leq\lambda}\dim
 A(\tilde\lambda,M),\end{equation} where the sum is certainly finite.

We define the $\lambda$-oldspace to be
$$A_{\op{old}}(\lambda,M):=\spn \{f(dz)|f\in A(\lambda,K)\quad
Kd|M\quad K\neq M \}.$$ This is contained in $A(\lambda,M)$ by the
$\SL (\R)$-invariance of $\Delta_\Gamma$, and the fact that $f(dz)$ is
$\Gamma_{\!0}(M)$-invariant when $f(z)$ is
$\Gamma_{\!0}(K)$-invariant and $Kd|M$.  We then define the
$\lambda$-newspace to be the orthogonal complement in $A(\lambda,M)$
with respect to the inner product$$(f,g)=\int_{F_{\Gamma_{\!0}(M)}}f(z)\overline{g(z)}d\mu(z),$$
i.e. $$A_{\op{new}}(\lambda,M):=A(\lambda,M)\ominus
A_{\op{old}}(\lambda,M).$$ 

We then define new spectral counting functions
\begin{eqnarray*}
N^{\op{old}}_{\Gamma_{\!0}(M)}(\lambda)& := & 1+\sum_{0<\tilde\lambda\leq\lambda}\dim
 A_{\op{old}}(\tilde\lambda,M)\qquad M>0\\
N^{\op{new}}_{\Gamma_{\!0}(M)}(\lambda)& := & \sum_{0<\tilde\lambda\leq\lambda}\dim
 A_{\op{new}}(\tilde\lambda,M)\qquad M>0.\\
\end{eqnarray*}For $M=1$ we of course define
 $N^{\op{old}}_{\Gamma_{\!0}(1)}(\lambda)=0$ and
 $N^{\op{new}}_{\Gamma_{\!0}(1)}(\lambda)=N_{\Gamma_{\!0}(1)}(\lambda)$.

\section{Calculation of asymptotic densities}
In order to calculate the main terms of $N^{\op{new}}_{\Gamma_{\!0}(M)}$ we
remind about some well known structure theory of arithmetical
functions.
When $f,g:\N\to\C$  are arithmetical functions we
define the \emph{Dirichlet convolution}, $f*g:\N\to\C$ to be the
arithmetical function
$$(f*g)(n)=\sum_{d|n}f(d)g\left(\frac{n}{d}\right).$$ We say that $f$
is multiplicative if $f(mn)=f(m)f(n)$ whenever  $(m,n)=1$. The
structure theory we shall use is the following: 
\begin{Theorem}\label{foldning} The arithmetical functions form a commutative group
  under Dirichlet convolution. The identity element is the
  function$$\begin{array}{rccc}
I : & \N    & \rightarrow  & \C\\
         & n & \mapsto      &
  \left[\frac{1}{n}\right]=\left\{\begin{array}{ll}1& \textrm{if }n=1\\ 0 & \textrm{otherwise.}\end{array}\right.
\end{array}$$
The multiplicative arithmetical functions form a subgroup.
\end{Theorem}
\begin{proof}
This follows from \cite[Theorems 2.6,2.8,2.14, 2.16]{MR55:7892}
\end{proof}
\begin{Example}(See \cite[\S 2.13]{MR55:7892} for details.)
Consider the arithmetical function $$\sigma_\alpha(n)=\sum_{d|n}d^\alpha. $$ Then
this in a \emph{multiplicative} arithmetical function whose inverse may be
calculated to be 

\begin{equation}\label{invers}\sigma_\alpha^{-1}(n)=\sum_{d|n}d^{\alpha}\mu(d)\mu\left(\frac{n}{d}\right),\end{equation}where
$\mu$ is the M\"{o}bius function, i.e.
$$\mu(n)=\left\{\begin{array}{ll}1 & \textrm{if }n=1\\(-1)^k & \textrm{if
    }n= p_1\cdots p_k \\0&\textrm{otherwise.}\end{array}\right.$$
The Mangoldt $\Lambda$-function
$$\Lambda(n)=\left\{\begin{array}{ll}\log p & \textrm{if
    }n=p^m\textrm{where $p$ is a prime and $m\geq1$}\\0&\textrm{otherwise.}\end{array}\right.$$
is an example of a non-multiplicative function. Another multiplicative
arithmetical function we will use is Eulers totient function
$$\Phi(n)=\#\{d\in\N|1\leq d\leq n \wedge (d,n)=1\}.$$

\end{Example}

We can now begin to calculate asymptotic densities of newforms.
We cite a result from \cite{andreas}.

\begin{Lemma}\label{dimensioner}$$\dim
  A(\lambda,\cdot)=\sigma_0*\dim A_{\op{new}}(\lambda,\cdot).$$\end{Lemma}
\begin{proof}This is Theorem 4.6.c) in Chapter III of \cite{andreas}.\end{proof}

Let now $f_i$, $i=1\ldots n$  be real positive functions of
decreasing order i.e $$f_{i+1}=o(f_{i})\textrm{ for }i=1\ldots n-1.$$ 

\begin{Prop}\label{rekurrens} Assume that for any $M\in \N$
  $$N_{\Gamma_{\!0}(M)}(\lambda)=\sum_{i=1}^{n-1}c_i(M)f_i(\lambda)\quad+O(f_n(\lambda)).$$ 
Then 
  $$N^{\op{new}}_{\Gamma_{\!0}(M)}(\lambda)=\sum_{i=1}^{n-1}c_i^{\op{new}}(M)f_i(\lambda)\quad+O(f_n(\lambda)).$$ 

where $c^{\op{new}}_i=c_i*\sigma_0^{-1}$.

\end{Prop}
\begin{proof}The $M=1$ case is clear by the definitions of
  $N^{\op{new}}_{\Gamma_{\!0}(1)}(\lambda)$ and $c_i^{\op{new}}(1)$. We observe
  that by lemma \ref{dimensioner} we have

\begin{eqnarray*}N_{\Gamma_{\!0}(M)}(\lambda)&=&1+\sum_{K|M}\sigma_0\left(\frac{M}{K}\right)\sum_{0<\tilde\lambda\leq\lambda}\dim
 A_{\op{new}}(\tilde\lambda,K)\\
&=&\sum_{K|M}\sigma_0\left(\frac{M}{K}\right)N_{\Gamma_{\!0}(K)}^{\op{new}}(\lambda).
\end{eqnarray*}
By the definition of $c_i^{\op{new}}$ we have 
$$c_i(M)=\sum_{K|M}\sigma_0\left(\frac{M}{K}\right)c_i^{\op{new}}(K)$$
 and therefore
\begin{align*}
\left|N^{\op{new}}_{\Gamma_{\!0}(M)}(\lambda)-\sum_{i=1}^{n-1}c^{\op{new}}_i(M)f_i(\lambda)\right|\leq&\left|N_{\Gamma_{\!0}(M)}(\lambda)-\sum_{i=1}^{n-1}c_i(M)f_i(\lambda)\right|\\
+\sum_{K|M\atop K\neq M}\sigma_0\left(\frac{M}{K}\right)&\left|N^{\op{new}}_{\Gamma_{\!0}(K)}(\lambda)-\sum_{i=1}^{n-1}c_i^{\op{new}}(K)f_i(\lambda)\right|.
\end{align*}
Induction in $M$ now gives that this is $\leq Cf_n(\lambda)$ which is
the desired result.
\end{proof}

The above proposition together with Theorem \ref{foldning} shows that
$c^{\op{new}}_i(N)$ is multiplicative if and only if $c_i(N)$ is
multiplicative. It also shows that if we know the expansion of the
counting function for eigenvalues of $\Delta_{\Gamma_{\!0}(M)}$ for
any $M\in\N$ and if these expansions are of the same type then it is easy to find the expansion of the
corresponding counting function for newforms. Finding the expansion of
$N_{\Gamma_{\!0}(M)}$ is the objective of the next section. 

\section{The refined Weyl law}
We start by citing a result by \citeasnoun{MR85j:11060b}(Theorem 5.2.1) which is the
basis of our calculations. We have adjusted the theorem to our
situation, and corrected the obvious misprint in the $O$-term.

\begin{Theorem}\label{venkov} The following asymptotic formula holds:
\begin{align*}N_{\Gamma_{\!0}(M)}(\lambda)&-\frac{1}{4\pi}\int_{-T}^T{\frac{\phi_M'}{\phi_M}\left(\frac{1}{2}+ir\right)}dr\\=&\frac{|F_M|}{4\pi}\lambda
  -\frac{k_M}{\pi}\sqrt{\lambda}\ln \sqrt{\lambda}+\frac{k_M(1-\ln 2)}{\pi}\sqrt{\lambda} 
  +  O(\sqrt{\lambda}/ \ln{\sqrt{\lambda}} )
\end{align*}
where $\phi_M$ is the determinant of the scattering matrix, $\lambda=1/4+T^2$, $k_M$ is the number of cusps of
$\Gamma_{\!0}(M)$ and $|F_M|$ is the area of the fundamental domain of $\Gamma_{\!0}(M)$.
\end{Theorem}
For definitions of cusps and the scattering matrix we refer to \cite{MR96f:11078} or \cite{MR55:2759}.

In order to apply Proposition \ref{rekurrens} to the estimate obtained
in Theorem \ref{venkov} we need to estimate the integral and therefore
also the logarithmic derivative
of the scattering matrix. We can do this by using the following
theorem which was proved by \citeasnoun{MR87e:11072}.

\begin{Theorem}
  Let $\phi_M(s)$ be the determinant of the scattering matrix for the Hecke congruence
  group of level M, $\Gamma_{\!0}(M)$, and let $\Lambda_\chi$ be the
  completed $L$-function of an  Dirichlet character mod $K$ ,$\chi$ ,
  i.e.
  $$\Lambda_\chi(s)=\Gamma\left(\frac{s}{2}\right)\sum_{n=1}^{\infty}\frac{\chi(n)}{n^s}\textrm{ \quad when }\Re (s)>1.$$ Then 
$$\phi_M(s)=(-1)^l\left(\frac{A(M)}{\pi^{k_M}}\right)^{1-2s}\prod_{i=1}^{k_M}\frac{\Lambda(2-2s,\overline{\chi_i})}{\Lambda(2s,\chi_i)}$$
where $l\in \N$, the $\chi_i$'s are some Dirichlet characters mod $K$ where
$K|M$, and $$A(M)=\prod_{\chi\op{primitive}\mod q\atop q|m,\, mq|M}\frac{qM}{(m,M/m)}.$$ 
The set $\{\chi_i|i=1,\ldots,k_M\}$ is closed under complex conjugation.
\end{Theorem}

We now use this to evaluate  the integral in theorem
\ref{venkov}. We let $B(M)=\frac{A(M)}{\pi^{k_M}}$. From the above we conclude that 
$$\frac{\phi_M'}{\phi_M}\left(\frac{1}{2}+ir\right)=-2\left(\ln B(M) +
  \sum_{i=1}^{k_M}\frac{\Lambda_{\chi_i}'}{\Lambda_{\chi_i}}(1-2it)+\frac{\Lambda_{\chi_i}'}{\Lambda_{\chi_i}}(1+2it)\right).$$ An easy consideration then shows that
$$
-\frac{1}{4\pi}\int_{-T}^T{\frac{\phi_M'}{\phi_M}\left(\frac{1}{2}+ir\right)}dr=\frac{T}{\pi}\ln
B(M)+\sum_{i=1}^{k_M}\frac{1}{\pi}\int_{-T}^T\frac{\Lambda_{\chi_i}'}{\Lambda_{\chi_i}}(1+2ir)dr.$$
We must therefore evaluate
$$\int_{-T}^T\frac{\Lambda_{\chi_i}'}{\Lambda_{\chi_i}}(1+2ir)dr,$$
and we  observe that 
$$\int_{-T}^T\frac{\Lambda_{\chi_i}'}{\Lambda_{\chi_i}}(1+2ir)d=\frac{1}{2}\int_{-T}^T\frac{\Gamma'}{\Gamma}\left(\frac{1}{2}+ir\right)dr+\int_{-T}^T\frac{L_\chi'}{L_\chi}\left(1+i2r\right)dr.$$
We shall address each term separately. To evaluate the first term we
use Stirling's approximation formula i.e. 
$$\frac{\Gamma'}{\Gamma}(s)=\log(s)-\frac{1}{2s}+O(|s|^{-2}),$$ valid
for $|\arg(s)-\pi|>\epsilon$. We see that for $|r|>\epsilon$ we have 
\begin{align*}
&\left|\frac{\Gamma'}{\Gamma}\left(\frac{1}{2}+ir\right)-\left(\log|r| +i \arg\left(\frac{1}{2}+ir\right)-\left(1+i2r\right)^{-1}\right)\right|\\
&\qquad\qquad\qquad\leq\left|\frac{\Gamma'}{\Gamma}\left(\frac{1}{2}+ir\right)-\left(\log\left|\frac{1}{2}+ir\right| +i
    \arg\left(\frac{1}{2}+ir\right)-\left(1+i2r\right)^{-1}\right)\right|\\
&\qquad\qquad\qquad+\left| \log\left|\frac{1}{2}+ir\right| -
  \log|r|\right|.
\end{align*}It is easy to see check that the last summand is $O((|r|\log
|r|)^{-1})$ while the first is $O(|r|^{-2})$ by Stirling's
 approximation formula. Hence
\begin{align*}
\frac{1}{2}\int_{-T}^{T}\frac{\Gamma'}{\Gamma}&\left(\frac{1}{2}+ir\right)dr\\
=&\frac{1}{2}\int_{-T\atop|r|>\epsilon}^{T}\log|r|
+i \arg\left(\frac{1}{2}+ir\right)-\left(1+i2r\right)^{-1}dr+\quad
O\left(\int_{\epsilon}^{T}\frac{1}{r\log r}\right)
\end{align*}The integral over $(1+i2r)^{-1}$ is bounded and the
integral over $i\arg(1/2+ir)$ vanishes. We conclude that
$$\frac{1}{2}\int_{-T}^{T}\frac{\Gamma'}{\Gamma}\left(\frac{1}{2}+ir\right)dr=T\log
T-T + O(\log(\log T))$$

To evaluate the integral over the logarithmic derivative of
$L_\chi(1+2ir)$ we note that
$$\int_\epsilon^T\frac{L_\chi'}{L_\chi}(1+2ir)dr=-i(\log
L_\chi(1+2iT)) +C$$ where $C$ is a constant and that the first term is
$O(\log T)$ by \cite[Theorem 12.24]{MR55:7892}. We conclude that 
$$-\frac{1}{4\pi}\int_{-T}^{T}\frac{\phi_M'}{\phi_M}\left(\frac{1}{2}+ir\right)dr=\frac{T}{\pi}\log
B(M)+\frac{k_M}{\pi}(T\log
T-T)+O_M(\log(T))$$

We have hence proven the following
\begin{Theorem}\label{refined}
The counting function $N_{\Gamma_{\!0}(M)}(\lambda)$ satisfies the
following asymptotic formula
\begin{align*}N_{\Gamma_{\!0}(M)}(\lambda)=&\frac{|F_M|}{4\pi}\lambda-\frac{2k_M}{\pi}\sqrt{\lambda}\log{\sqrt{\lambda}}\\&+\frac{1}{\pi}\left[(2-\log
2+\log\pi)k_M)-\log(A(M))\right]\sqrt{\lambda}+O(\sqrt{\lambda}/\log{\sqrt{\lambda}})
\end{align*}
\end{Theorem}
This theorem puts us in a situation where proposition \ref{rekurrens}
can be applied with 
\begin{eqnarray*}f_1(\lambda)&=&\lambda\\
f_2(\lambda)&=&\sqrt{\lambda}\log{\sqrt{\lambda}}\\
f_3(\lambda)&=&\sqrt{\lambda}\\
f_4(\lambda)&=&\sqrt{\lambda}/\log{\sqrt{\lambda}}.
\end{eqnarray*}

From \cite[Theorem 1.43]{MR47:3318} we conclude that \begin{align}
\label{antalcusps} k_M&=\sum_{d|M}\Phi((d,M/d))\\
\label{areal} |F_M|&=\frac{\pi}{3}M\prod_{p|M\atop p\, \op{prime}}(1+p^{-1}).\end{align}
This means that we have explicit expressions for all the terms in
theorem \ref{refined} except $A(M)$. We need to know the number of
primitive Dirichlet characters mod $K$. We hence define
$$D(K)=\#\{\chi \textrm{ primitive Dirichlet character mod }K\}.$$
Then we have
\begin{Lemma}The arithmetical function $D(K)$ is multiplicative and
  satisfies
$$D(K)=(\Phi*\mu)(K)$$
\end{Lemma}
\begin{proof}
From \cite{MR55:7892} theorem 6.15 and theorem 8.18 we conclude that
$\Phi(K)=\sum_{d|K}D(d)=(u*D)(K)$ where $u(n)=1$ for $n\in \N$. Since
$\Phi$ and $u$ are multiplicative we use theorem \ref{foldning} to
conclude that $D$ is multiplicative. Theorem 2.1 in \cite{MR55:7892}
proves that $u^{-1}=\mu$ so $$\Phi*\mu=u*D*\mu=u*u^{-1}*D=D$$ which
concludes the proof 
\end{proof}
\section{Coefficients related to newforms}
In this section we calculate  $c_1^{\op{new}}$,$c_2^{\op{new}}$ and  $c_3^{\op{new}}$. The calculations
are basically corollaries of proposition \ref{rekurrens} and theorem
\ref{refined}. We are particularly interested in the case when $c_2^{\op{new}}(M)=c_3^{\op{new}}(M)=0$. If this is the case we
  say that $N^{\op{new}}_{\Gamma_{\!0}(M)}$ 
 \emph{is of cocompact type}. To see why this is sensible we remind of the following

\begin{Theorem}\label{compactcase}Assume $\Gamma$ is a cocompact Fuchsian group and let
  $N_\Gamma(\lambda)$ be the counting function for the eigenvalues of
  $\Delta_\Gamma$. Then   $$N_\Gamma(\lambda)=\frac{|F_\Gamma|}{4\pi}\lambda+O(\sqrt{\lambda}/\log\sqrt{\lambda}),$$
  where $|F_\Gamma|$ is the area of an fundamental domain of $\Gamma$. \end{Theorem}
\begin{proof}
This is a special case of \cite{MR85j:11060b}(Theorem 5.2.1). Notice
again that we have corrected the obvious misprint. 
\end{proof}

Hence  $N^{\op{new}}_{\Gamma_{\!0}(M)}$ is of cocompact type if and
only if it 'has the same shape' as if it where the counting function
of the eigenvalues related to a cocompact group.

\subsection{The first coefficient}We start by calculating
$c_1^{\op{new}}(M)$. This is the simplest of the three coefficients.
\begin{Prop}\label{1koefficient}
The arithmetical function $v(M)=12c_1^{\op{new}}(M)$ is multiplicative
and satisfies\begin{equation}\label{v(n)}v(p^n)=\left\{\begin{array}{ll}
1&\textrm{ if }n=0\\
p-1&\textrm{ if }n=1\\
p^2-p-1&\textrm{ if }n=2\\
(p^3-p^2-p+1)p^{n-3}&\textrm{ if }n\geq 3\end{array}\right.\end{equation} when
$p$ is a prime. We furthermore have    
$$L_v(s):=\sum_{n=1}^\infty\frac{v(n)}{n^s}=\frac{\zeta(s-1)}{\zeta(2s)\zeta(s)},$$
where $\zeta(s)$ is Riemann's zeta function.
\end{Prop}
\begin{proof}By using proposition \ref{rekurrens}, theorem
  \ref{refined} and \ref{areal} we conclude that $$ M\prod_{p|M\atop
  p\, \op{prime}}(1+p^{-1})=(\sigma_0*v)(M).$$ Since the left hand
  side and $\sigma_0$ are multiplicative theorem \ref{foldning} says
  that $v$ is multiplicative. By considering the case where $M=p^m$ we
  see that
  $$p^m+p^{m-1}=\sum_{d|p^m}\sigma_0(d)v\left(\frac{p^m}{d}\right)=\sum_{i=0}^{n}(i+1)v(p^{n-i}).$$ By applying the theory of generating functions to this relation we find that if $$f_p(c)=\sum_{n=0}^\infty v(p^n)x^n\quad\textrm{ then }\quad f_p(x)=\frac{(1-x^2)(1-x)}{1-px}.$$ By making formal expansion we get (\ref{v(n)}). Since $v$ is multiplicative the claim about $L_v$ follows. 
\end{proof}

\subsection{The second coefficient} We now calculate
$c_2^{\op{new}}(M)$. We remind that by proposition \ref{rekurrens} and theorem
\ref{refined} we have $$c_2^{\op{new}}(M)=-\frac{2}{\pi}(k_{(\cdot)}*\sigma_0^{-1})(M)$$
We hence need to have more information about the number of cusps of
$\Gamma_{\!0}(M)$ 
\begin{Lemma}\label{k-egenskaber}
  The number of cusps, $k_M$, of $\Gamma_{\!0}(M)$ is a multiplicative
  arithmetical function and satisfies
\begin{equation}k_{p^m}\left\{\begin{array}{ll}
1&\textrm{ if }m=0\\
2&\textrm{ if }m=1\\
2p^{n}&\textrm{ if }m=2n+1\textrm{ where }n>1\\
(p+1)p^{n-1}&\textrm{ if }m=2n\textrm{ where }n>1 \end{array}\right.\end{equation}

\end{Lemma}
\begin{proof}
We noted earlier in (\ref{antalcusps}) that $$k_M=\sum_{d|M}\Phi((d,M/d)).$$ Let $M_1,M_2\in
\N$ and assume $(M_1,M_2)=1$. Then 
\begin{align*}
  k_{M_1M_2}=&\sum_{d|M_1M_2}\Phi((d,(M_1M_2)/d))\\
=&\sum_{d_1|M_1}\sum_{d_2|M_2}\Phi((d_1d_2,(M_1M_2)/(d_1d_2)))\\
=&\sum_{d_1|M_1}\sum_{d_2|M_2}\Phi((d_1,M_1/d_1)(d_2,M_2/d_2))\\
=&\sum_{d_1|M_1}\Phi((d_1,M_1/d_1))\sum_{d_2|M_2}\Phi((d_2,M_2/d_2))\\
=&k_{M_1}k_{M_2}
\end{align*} Hence $k_{M}$ is multiplicative. The claim about
$k_{p^m}$ is clear for $m=0$ and $m=1$. Assume $m\geq 2$. We then have
\begin{align*}
k_{p^m}&=\sum_{i=0}^m\Phi((p^i,p^{m-i}))\\
       &=\sum_{i=0}^m\Phi(p^{\min(i,m-i)})\\
       &=2+\sum_{i=1}^{m-1}(p-1)p^{\min(i,m-i)-1}\\
\intertext{We now assume $m=2n+1$.}
       &=2+(p-1)\left(\sum_{i=1}^{n}p^{i-1}+\sum_{i=n+1}^{2n}p^{2n-i}\right)\\
       &=2+2(p-1)\sum_{i=0}^{n-1}p^i\\
       &=2+2(p-1)\frac{1-p^n}{1-p}=2p^n.
\end{align*}The even case is similar.
\end{proof}
From the above we can now prove the following
\begin{Prop}
  The second coefficients ,$c_2^{\op{new}}(M)$, is a multiplicative
  arithmetical function and satisfies
\begin{equation}-\frac{\pi}{2}c_2^{\op{new}}(p^m)=\left\{\begin{array}{ll}
1&\textrm{ if }m=0\\
0&\textrm{ if }m=2n+1\\
p-2&\textrm{ if }m=2\\
(p+1)^2p^{n-1}&\textrm{ if }m=2n\textrm{ where }n>1 \end{array}\right.\end{equation}\end{Prop}
\begin{proof}
  From lemma \ref{k-egenskaber} and theorem \ref{foldning} follows
  that $c_2^{\op{new}}(M)$ is  multiplicative.  From (\ref{invers}) it
  is easy to see that $$\sigma_0^{-1}(p^m)=\left\{\begin{array}{rl}
1&\textrm{ if }m=0\\
-2&\textrm{ if }m=1\\
1&\textrm{ if }m=2\\
0&\textrm{ otherwise. }\end{array}\right.$$ Hence 
$$c_2^{\op{new}}(p^m)=-\frac{2}{\pi}(k_{p^{m}}-2k_{p^{m-1}}+k_{p^{m}}),\textrm{
when }m\geq 2.$$ Using lemma \ref{k-egenskaber} it is now easy to
check the claim. We omit the details.
\end{proof}
As an easy corollary we get the following
\begin{Corollary}
  The second coefficient, $c_2^{\op{new}}(M)$, is non-zero  if and only if $M=t^2$ where $t\in \N$ is not
  of the form $t=2t'$ with $(2,t')=1$.
\end{Corollary}

\subsection{The third coefficient} We finally calculate
  $c_3^{\op{new}}(M)$. This is the most difficult of the three coefficients.

We start by observing that by proposition \ref{rekurrens} and theorem
\ref{refined} 

$$c_3^{\op{new}}(M)=\frac{1}{\pi}\left(\left(2-\log
    2+\log\pi\right)\left(-\frac{\pi}{2}c_2^{\op{new}}(M)\right)-L(M)\right)$$ where $$L(M)=\left(\log A(\cdot)*\sigma_0^{-1}\right)(M).$$ We hence direct our attention to $L(M)$.
\begin{Lemma}
  Assume $(M_1,M_2)=1$. Then $$L(M_1M_2)=U(M_1)L(M_2)+U(M_2)L(M_1)$$
  where $$U(M)=\sum_{d|M}\sum_{m|d}\sum_{q|(m,\frac{d}{m})}D(q)\sigma_0^{-1}\left(\frac{M}{d}\right).$$
\end{Lemma}
\begin{proof}We have 
\begin{align*}
&L(M_1M_2)=\sum_{d|M_1M_2}\log
         A(d)\sigma_0^{-1}\left(\frac{M_1M_2}{d}\right)\allowdisplaybreaks
         \\ 
         &=\sum_{d|M_1M_2}\sum_{q|m\atop
         mq|d}D(q)\log\left(\frac{qd}{(m,\frac{d}{m})}\right)\sigma_0^{-1}\left(\frac{M_1M_2}{d}\right)\allowdisplaybreaks
         \\
         &=\sum_{d|M_1M_2}\sum_{m|d}\sum_{q|(m,d/m)}D(q)\log\left(\frac{qd}{(m,\frac{d}{m})}\right)\sigma_0^{-1}\left(\frac{M_1M_2}{d}\right)\allowdisplaybreaks
         \\
         &=\sum_{d_1|M_1}\sum_{d_2|M_2}\sum_{m_1|d_1}\sum_{m_2|d_2}\sum_{q_1|(m_1,\frac{d_1}{m_1})}\sum_{q_2|(m_2,\frac{d_2}{m_2})}D(q_1q_2)\log\left(\frac{q_1q_2d_1d_2}{(m_1m_2,\frac{d_1d_2}{m_1m_2})}\right)\sigma_0^{-1}\left(\frac{M_1M_2}{d_1d_2}\right)\end{align*}
The summand is clearly 
$$D(q_1)D(q_2)\sigma_0^{-1}\left(\frac{M_1}{d_1}\right)\sigma_0^{-1}\left(\frac{M_2}{d_2}\right)\left( \log
  \left(\frac{q_1d_1}{(m_1,\frac{d_1}{m_1})}\right)+\log\left(\frac{q_2d_2}{(m_2,\frac{d_2}{m_2})}\right)\right).$$
We have 
\begin{align*}
\sum_{d_1|M_1}&\sum_{d_2|M_2}\sum_{m_1|d_1}\sum_{m_2|d_2}\sum_{q_1|(m_1,\frac{d_1}{m_1})}\sum_{q_2|(m_2,\frac{d_2}{m_2})}\!\!\!\!\!\!D(q_1)D(q_2)\sigma_0^{-1}\!\!\left(\frac{M_1}{d_1}\right)\!\sigma_0^{-1}\!\!\left(\frac{M_2}{d_2}\right)\!\!\left(\! \log
 \! \left(\frac{q_1d_1}{(m_1,\frac{d_1}{m_1})}\right)\!\right)\\
&=U(M_2)L(M_1),
\end{align*}
from which the identity easily follows. 
\end{proof}
It turns out that $U$ is a very nice arithmetical function. In fact we
have the following.
\begin{Lemma}\label{U-funktionen}
   The function $U(M)$, is a multiplicative
  arithmetical function and satisfies
\begin{equation}U(p^m)=\left\{\begin{array}{ll}
1&\textrm{ if }m=0\\
0&\textrm{ if }m=2n+1\\
p-2&\textrm{ if }m=2\\
(p^2-2p+1)p^{n-2}&\textrm{ if }m=2n\textrm{ where }n>1 \end{array}\right.\end{equation}
\end{Lemma}

\begin{proof}
  Let $M_1, M_2\in \N$ be coprime. Then 
  \begin{align*}
    U&(M_1M_2)=\sum_{d|M_1M_2}\sum_{m|d}\sum_{q|(m,\frac{d}{m})}D(q)\sigma_0^{-1}\left(\frac{M_1M_2}{d}\right)\allowdisplaybreaks \\
&=\sum_{d_1|M_1}\sum_{d_2|M_2}\sum_{m_1|d_1}\sum_{m_2|d_2}\sum_{q_1|(m_1,\frac{d_1}{m_1})}\sum_{q_2|(m_2,\frac{d_2}{m_2})}D(q_1)D(q_2)\sigma_0^{-1}\left(\frac{M_1}{d_1}\right)\sigma_0^{-1}\left(\frac{M_2}{d_2}\right)\allowdisplaybreaks
\\
&=U(M_1)U(M_2).
  \end{align*} Hence U is multiplicative.

Let $p$ be a prime and $m\in \N$. We assume $m\geq 2$ Then
\begin{align*}
  L(p^m)&=\sum_{i=0}^m\sum_{j=0}^i\sum_{l=0}^{\min(j,i-j)}D(p^l)\sigma_0^{-1}\left(p^{m-i}\right)\allowdisplaybreaks\\
=&\sum_{j=0}^{m-2}\sum_{l=0}^{\min(j,m-2-j)}D(p^l)-2\sum_{j=0}^{m-1}\sum_{l=0}^{\min(j,m-1-j)}D(p^l)+\sum_{j=0}^m\sum_{l=0}^{\min(j,m-j)}D(p^l)\allowdisplaybreaks\\
\intertext{Assume $j\leq n-2-j$. Then $j\leq n-1-j\leq n-j$ and we have that all
  minimum values are $j$. Hence these terms cancels out. We now assume
$m=2n+1$. Hence we may sum from $j\geq(2n+1)/2-1=n-1/2$.}
=&\sum_{j=n}^{m-2}\sum_{l=0}^{\min(j,m-2-j)}D(p^l)-2\sum_{j=n}^{m-1}\sum_{l=0}^{\min(j,m-1-j)}D(p^l)+\sum_{j=n}^m\sum_{l=0}^{\min(j,m-j)}D(p^l)\allowdisplaybreaks\\
=&\sum_{j=n}^{m-2}\sum_{l=0}^{m-2-j}D(p^l)-2\sum_{j=n}^{m-1}\sum_{l=0}^{m-1-j}D(p^l)+\sum_{j=n+1}^m\sum_{l=0}^{m-j}D(p^l)+\sum_{l=0}^mD(p^l)\allowdisplaybreaks\\
=&\sum_{l=0}^{m-2-n}D(p^l)-2\sum_{l=0}^{m-1-n}D(p^l)+\sum_{l=0}^nD(p^l)\allowdisplaybreaks\\
 &+\sum_{j=n+1}^{m-2}\left(\sum_{l=0}^{m-2-j}D(p^l)-2\sum_{l=0}^{m-1-j}D(p^l)+\sum_{l=0}^{m-j}D(p^l)\right)\allowdisplaybreaks\\
 &-2\sum_{l=0}^{m-1-(m-1)}D(p^l)+\sum_{l=0}^{m-(m-1)}D(p^l)+\sum_{l=0}^{m-m}D(p^l)\allowdisplaybreaks\\
=&-D(p^n)+\sum_{j=n+1}^{m-2}(-2D(p^{m-1-j})+D(p^{m-1-j})+D(p^{m-j}))+D(p)\allowdisplaybreaks\\
=&-D(p^n)-\sum_{j=1}^{n-1}D(p^j)+\sum_{j=2}^{n}D(p^j)+D(p)=0
\end{align*}
The even case is similar but slightly easier. The $m=1$ case is also
similar.
\end{proof}
\begin{Remark}
  By successive use of the two lemmas above we find that $$L(p_1^{n_1}\ldots
p_k^{n_k})=\sum_{i=1}^k\left(\prod_{j\in\{1,\ldots,k\}\backslash \{ i \} }U(p_j^{n^j})\right)L(p_i^{n_i}),$$ when $p_1,\ldots,p_k$
are different primes. Notice that $L(p_i^{n_i})$ is of the form
$\tilde{m}_i\log p_i$ where $\tilde{m}_i\in\Z$. We also note that
$U(M)\in \Z$. Hence $L$ is on the form $$m_1\log p_1+\ldots,m_k\log
p_k\textrm{ where }m_i\in\Z.$$ By unique factorization in $\N$ this is
zero if and only if $m_i=0$ for all $i$'s. We would therefore like to know
when $L(p^m)$ is zero. 
\end{Remark}

\begin{Lemma}
  The function $L(p^m)$ satisfies
\begin{equation}L(p^m)=\left\{\begin{array}{ll}
2\left(\sum_{j=0}^nD(p^j)\right)\log p&\textrm{ if }m=2n+1\\
\left(\sum_{j=0}^{n-1}D(p^j)+mD(p^n)\right)\log p&\textrm{ if }m=2n.\end{array}\right.\end{equation}
In particular $L(p^m)$ is never zero.
\end{Lemma}
\begin{proof}
  This follows by a lengthy but elementary calculation similar to that
  in the proof of lemma \ref{U-funktionen}.
\end{proof}

From the above lemma and the preceding remark we conclude that $$L(p_1^{n_1}\ldots p_k^{n_k})=0$$
if and only if $U(p_i^{n_i})=0$ for at least two different
primes. Since $c_2^{\op{new}}(M)=c_3^{\op{new}}(M)=0$ if and only if
$c_2^{\op{new}}(M)=L(M)=0$ we have proved the following which settles
the question of when $N^{\op{new}}_{\Gamma_{\!0}(M)}(\lambda)$ is of
cocompact type.

\begin{Theorem}\label{cocompactform}
  Let $M\in\N$ and let $n,t\in\N$ be the integers defined uniquely by
  the requirements that $n$ should be squarefree and $M=t^2n$. Then
  $N^{\op{new}}_{\Gamma_{\!0}(M)}(\lambda)$ is of cocompact type if
  and only if $n,t$ satisfies one of the following:
  \begin{enumerate}
  \item{$n$ contains more than one prime.}
  \item{$n>1$ and $4|M$ and $(2,M/4)=1$.}
\end{enumerate}
\end{Theorem}

\section{Concluding remarks}
\begin{Remark}
From proposition \ref{1koefficient} we conclude that 
$$N_{\Gamma_{\!0}(M)}^{\op{new}}(\lambda)=\frac{1}{12}\lambda+O(\sqrt{\lambda}\log\sqrt{\lambda})$$
  if and only if $M\in\{1,2,4\}$. This shows that theorem 2 of \cite{MR99g:11069}
  cannot be generalized to more general Hecke
  congruence groups  by simply choosing another character.
\end{Remark}
\begin{Remark}
  We wish to draw attention to a particular case of theorem
  \ref{cocompactform} namely the case when $M>1$ is
  squarefree with an even number of primes. 
Hence, by Theorem \ref{cocompactform} (1) $N_{\Gamma_{\!0}(N)}^{\op{new}}(\lambda)$ has the same
form as if it were the counting function
for the eigenvalues related to a co-compact group with invariant area
$4\pi c_1^{\op{new}}(M)$. We can give an alternative and much more
  sophisticated proof of this by referring to the
Jacquet-Langlands correspondence. A part of this correspondence is
described classically in \cite{MR1832538} where the following is proven:

\begin{quote} Let $\mathcal{O}$ be a maximal order in an indefinite rational
quaternion division algebra over $\Q$, and let $d=d(\mathcal{O})$ be
its (reduced) discriminant. This is always a \emph{squarefree integer with
an even number of prime factors}. The norm one unit group
$\mathcal{O}^1$ can be viewed as a Fuchsian group \emph{which is
  cocompact}. Then:\end{quote}

\begin{quote}\emph{The eigenvalues of the Laplacian on
  $\mathcal{O}^1\setminus\mathcal{H}$ are exactly the same (with
  multiplicities) as the eigenvalues corresponding to the newspace on
  $\Gamma_{\!0}(d)\setminus\mathcal{H}$.}
\end{quote}

Hence $N_{\Gamma_{\!0}(N)}^{\op{new}}(\lambda)$ \emph{is} the counting
function for the eigenvalues related to a cocompact group, and hence
obviously has the corresponding type as predicted by Theorem \ref{compactcase}.
This has our theorem \ref{cocompactform} as an easy corollary (compare
with theorem \ref{compactcase}). We note that any squarefree $d$ with an even number of primes may be
constructed in this way. 

Our calculation indicates that there might be a similar correspondence in a lot of
other cases. We hope to address this on a later occasion.
\end{Remark}

\section*{Acknowledgment}I am grateful to Professor A.B. Venkov and
Professor E. Balslev for inspiring conversations and helpful advice. 

\harvardparenthesis{round}
\bibliographystyle{dcu}
\bibliography{bibliografi}

\end{document}